\documentclass[12pt]{article}
\title{Generic section of a hyperplane arrangement and twisted Hurewicz maps}
\author{Masahiko Yoshinaga\\
Department of Mathematice, Graduate School of Science, \\
Kobe University, 
1-1 Rokkodai, Kobe 657-8501 Japan\\
email: myoshina@math.kobe-u.ac.jp}
\date{\today}

\usepackage{amssymb}

\newtheorem{Def}{Definition}

\newtheorem{Thm}[Def]{Theorem}
\newtheorem{Lemma}[Def]{Lemma}
\newtheorem{Cor}[Def]{Corollary}
\newtheorem{Rem}[Def]{Remark}

\newcommand{\bbC}{\mathbb{C}}

\newcommand{\bbP}{\mathbb{P}}

\newcommand{\bbZ}{\mathbb{Z}}

\newcommand{\rmD}{\mathrm{D}}

\newcommand{\rmS}{\mathrm{S}}

\newcommand{\calA}{\mathcal{A}}

\newcommand{\calL}{\mathcal{L}}

\newcommand{\owari}{\hfill$\square$}

\makeatletter
\def\rank{\mathop{\operator@font rank}\nolimits}
\def\Hom{\mathop{\operator@font Hom}\nolimits}
\makeatother

\begin{document}
\maketitle

\begin{abstract}
We consider a twisted version of the Hurewicz map on the complement 
of a hyperplane arrangement. The purpose of this paper is to 
prove surjectivity of the twisted Hurewicz map under some 
genericity conditions. 
As a corollary, 
we also prove that a generic section of the complement of a 
hyperplane arrangement has non-trivial homotopy groups. 
\end{abstract}

\section{Twisted Hurewicz map}

Let $X$ be a topological space with a base point $x_0\in X$ 
and $\calL$ a local system of $\bbZ$-modules on $X$. 
Let $f:(S^n, *)\rightarrow (X, x_0)$ 
be a continuous map from the sphere $S^n$ with $n\geq 2$. 
Since $S^n$ is simply connected, the pullback 
$f^*\calL$ turns out to be a trivial 
local system. Thus given a local section $t\in\calL_{x_0}$, 
$f\otimes t$ determines a twisted cycle with coefficients in $\calL$. 
This induces a twisted version of the Hurewicz map: 
$$
h:\pi_n(X, x_0)\otimes_\bbZ\calL_{x_0}
\longrightarrow 
H_n(X, \calL). 
$$
The classical Hurewicz map is corresponding to the case of 
trivial local system $\calL=\bbZ$ with $t=1$. 

\section{Main result}

Let $\calA$ be an essential affine hyperplane arrangement 
in an affine space $V=\bbC^\ell$, with $\ell\geq 3$. 
Let $M(\calA)$ denote the complement 
$V -\bigcup\nolimits_{H\in\calA}H$. 
A hyperplane $U\subset V$ is said to be {\it generic to} $\calA$ if 
$U$ is transversal to the stratification induced from $\calA$. 
Let $i:U\cap M(\calA)\hookrightarrow M(\calA)$ denote the 
inclusion. 

In this notation, the main result of this paper is the following: 

\begin{Thm}
\label{thm:main}
Let $\calL':= i^*\calL$ be the 
restriction of a nonresonant local system 
$\calL$ of arbitrary rank on $M(\calA)$. 
Then the twisted Hurewicz map 
$$
h:\pi_{\ell-1}(U\cap M(\calA), x_0)\otimes_\bbZ\calL_{x_0}
\longrightarrow 
H_{\ell-1}(U\cap M(\calA), \calL') 
$$
is surjective. 
\end{Thm}

For the notion of ``nonresonant local system'', see Theorem \ref{thm:vanishing}.

Theorem \ref{thm:main} should be compared with a result proved by 
Randell in \cite{ran-hom}. He proved that the Hurewicz 
homomorphism $\pi_k(M(\calA))\rightarrow H_k(M(\calA), \bbZ)$ 
is equal to the zero map when $k\geq 2$ for any $\calA$. 
However little is known about twisted Hurewicz maps for 
other cases.

The key ingredient for our proof of Theorem \ref{thm:main} is 
an affine Lefschetz theorem of Hamm, 
which asserts that 
$M(\calA)$ has the homotopy type of a finite CW complex whose 
$(\ell-1)$-skeleton has the homotopy type of $U\cap M(\calA)$. 
We obtain $(\ell -1)$-dimensional spheres in $U\cap M(\calA)$ 
as boundaries of the $\ell$-dimensional top cells. 
Applying a vanishing theorem for local system homology groups, 
we show that these spheres generate 
the twisted homology group $H_{\ell -1}(U\cap M(\calA), \calL)$. 
We should note that the essentially same arguments are used 
in \cite{dp-eq} to compute the rank of 
$\pi_{\ell-1}(U\cap M(\calA), x_0)\otimes_\bbZ\calL_{x_0}$ under 
a certain asphericity condition on $\calA$. 


\section{Topology of complements}
\label{sec:top}

The cell decompositions of affine varieties or hypersurface 
complements are well studied subjects. 
Let $f\in\bbC [x_1, \cdots, x_\ell]$ 
be a polynomial and $D(f):=\{x\in\bbC^\ell |\ f(x)\neq 0\}$ be 
the hypersurface complement defined by $f$. 

\begin{Thm}
\label{thm:hamm}
\normalfont
({\bf Affine Lefschetz Theorem \cite{ham-lef}})
Let $U$ be a sufficiently generic hyperplane in $\bbC^\ell$. Then, 
\begin{itemize}
\item[(a)] The space $D(f)$ has the homotopy type of a space obtained 
from $D(f)\cap U$ by attaching $\ell$-dimensional cells. 
\item[(b)] Let 
$i_p:H_p(D(f)\cap U, \bbZ)\longrightarrow H_p(D(f), \bbZ)$ 
denote the homomorphism induced by the natural inclusion 
$i:D(f)\cap U\hookrightarrow D(f)$. Then 
$$
i_p \mbox{ is }
\left\{
\begin{array}{ll}
\mbox{isomorphic }&\mbox{for }p=0, 1, \cdots, \ell-2 \\
\mbox{surjective }&\mbox{for }p=\ell -1. 
\end{array}
\right.
$$
\end{itemize}
\end{Thm}

Suppose $i_{\ell-1}$ is also isomorphic. 
Then as noted by Dimca and Papadima \cite{dim-pap} (see 
also Randell \cite{ran-mor}), 
the number of 
$\ell$-dimensional cells attached would be equal to the Betti number 
$b_\ell(D(f))$ and the chain boundary map 
$\partial :C_\ell(D(f), \bbZ)\rightarrow C_{\ell-1}(D(f), \bbZ)$ 
of the cellular chain complex associated to the cell decomposition 
is equal to zero. Otherwise 
$i_{\ell-1}:H_{\ell-1}(D(f)\cap U, \bbZ)\longrightarrow 
H_{\ell-1}(D(f), \bbZ)$ 
has a nontrivial kernel $\partial(C_\ell(D(f), \bbZ))$. 

In the case of hyperplane arrangements, 
homology groups and homomorphisms $i_p$ are described combinatorially 
in terms of the intersection poset \cite{orl-ter1}. 
Let us recall some notation. 
Let $\calA$ be a finite set of affine hyperplanes in $\bbC^\ell$, 
$$
L(\calA)=\{ X=\bigcap_{H\in I}H\ |\ I\subset \calA\} 
$$
be the set of nonempty intersections of elements of $\calA$ 
with reverse inclusion $X<Y\Longleftrightarrow X\supset Y$, 
for $X, Y\in L(\calA)$. 
Define a rank function on $L(\calA)$ by
$$
r:L(\calA)\longrightarrow \bbZ_{\geq 0},\ X\longmapsto {\mathrm{codim}} X, 
$$
the M\"obius function $\mu:L(\calA)\longrightarrow \bbZ$ by 
$$
\mu(X)=
\left\{
\begin{array}{ll}
1 &\mbox{ for }X=V\\
-\sum_{Y<X}\mu(Y), &\mbox{ for } X>V, 
\end{array}
\right.
$$
and the characteristic polynomial $\chi(\calA, t)$ by 
$$
\chi(\calA, t)=\sum_{X\in L(\calA)}\mu(X)t^{\dim X}. 
$$
Let $E^1=\bigoplus_{H\in\calA}\bbC e_H$ and $E=\wedge E^1$ be 
the exterior algebra of $E^1$, with $p$-th graded term 
$E^p=\stackrel{p}{\wedge}E^1$. Define a $\bbC$-linear map 
$\partial :E\rightarrow E$ by $\partial 1=0$, $\partial e_H=1$ and 
for $p\geq 2$ 
$$
\partial(e_{H_1}\cdots e_{H_p})=\sum_{k=1}^p(-1)^{k-1}
e_{H_1}\cdots \widehat{e_{H_k}}\cdots e_{H_p}
$$
for all $H_1, \cdots, H_p\in\calA$. A subset $S\subset \calA$ is 
said to be dependent if $r(\cap S)<|S|$, where $\cap S=\cap_{H\in S}H$.  
For $S=\{H_1, \cdots, H_p\}$, we write 
$e_S:=e_{H_1}\cdots e_{H_p}$. 

\begin{Def}
\normalfont
\label{def:os}
Let $I(\calA)$ be the ideal of $E(\calA)$ generated by 
$$
\{e_S\ |\ \cap S=\phi\}\cup\{\partial e_S\ |\ S\mbox{ is dependent}\}. 
$$
The Orlik-Solomon algebra $A(\calA)$ is defined by 
$A(\calA)=E(\calA)/I(\calA)$. 
\end{Def}

\begin{Thm}
\normalfont
\label{thm:os}
({\bf Orlik-Solomon \cite{orl-sol}}) 
Fix a defining linear form $\alpha_H$ for each $H\in\calA$. Then 
the correspondence $e_H\mapsto d\log \alpha_H$ induces an isomorphism 
of graded algebras: 
$$
A(\calA)\stackrel{\cong}{\longrightarrow} H^*(M(\calA), \bbC). 
$$
The Betti numbers of $M(\calA)$ are given by 
$$
\chi(\calA, t)=\sum_{k=0}^{\ell}(-1)^k b_k(M(\calA))t^{\ell -k}. 
$$
\end{Thm}

From the above description of cohomology ring of $M(\calA)$, 
we have: 

\begin{Thm}
\normalfont
Let $\calA$ be a hyperplane arrangement in $\bbC^\ell$ and $U$ be 
a hyperplane generic to $\calA$. Then 
$i:U\cap M(\calA)\hookrightarrow M(\calA)$ induces isomorphisms 
$i_p:H_p(M(\calA)\cap U, \bbZ)\stackrel{\cong}{\longrightarrow} 
H_p(M(\calA), \bbZ)$ 
for $p=0, \cdots, \ell -1$. 
\end{Thm}

\noindent
{\bf Proof.} 
It is easily seen from the genericity that 
\begin{equation}
\label{eq:trunc}
L(\calA\cap U)\cong L_{\leq \ell -1}(\calA):=
\{ X\in L(\calA)|\ r(X)\leq \ell-1\}. 
\end{equation}
In particular a generic intersection preserves the part of 
rank $\leq \ell-1$. Hence $A(\calA\cap U)\cong A^{\leq \ell-1}(\calA)$. 
This induces isomorphisms 
$H^{\leq \ell-1}(M(\calA))\cong H^{\leq \ell-1}(M(\calA)\cap U)$. 
Since homology groups $H_*(M(\calA), \bbZ)$ are torsion free, 
the theorem is the dual of these isomorphisms. 
\owari

Using these results inductively, the complement $M(\calA)$ of 
the hyperplane arrangement $\calA$ has a minimal cell decomposition. 

\begin{Thm}
\normalfont
({\bf \cite{ran-mor}\cite{dim-pap}\cite{pap-suc}})
The complement $M(\calA)$ is homotopic to a minimal CW cell complex, 
i.e. the number of 
$k$-dimensional cells is equal to the Betti number 
$b_k(M(\calA))$ for each $k=0, \cdots, \ell$. 
\end{Thm}

\section{Proof of the main theorem}

First we recall the vanishing theorem of homology 
groups for a ``generic'' or nonresonant local system $\calL$ of 
complex rank $r$. 

Let $\calA$ be a hyperplane arrangement in $\bbC^\ell$, let 
$U$ be a hyperplane generic to $\calA$ and let 
$i:M(\calA)\cap U\hookrightarrow M(\calA)$ be the inclusion. 
Now we assume that 
$\calA$ is essential, i.e., 
$\calA$ contains linearly independent $\ell$ hyperplanes 
$H_1, \cdots, H_\ell \in\calA$. 

Let $\bbP^\ell$ be the projective space, which is a 
compactification of our vector space $V$. 
The projective closure of $\calA$ is defined as 
$\calA_\infty:=\{ \bar{H}|\ H\in\calA\}\cup \{H_\infty\}$, where 
$\bbP^\ell=V\cup H_\infty$. 
A non-empty intersection 
$X\in L(\calA_\infty)$ defines the subarrangement 
$(\calA_\infty)_X=\{H\in\calA_\infty|\ X\subset H\}$ of $\calA_\infty$. 
A subspace $X\in\calA_\infty$ is called dense if 
$(\calA_\infty)_X$ is indecomposable, that is, not the product of two 
non-empty arrangements. 
Let $\rho: \pi_1(M(\calA), x_0)\rightarrow GL_r(\bbC)$ be 
the monodromy representation associated to $\calL$. 
Choosing a point 
$p\in X\setminus\bigcup_{H\in\calA_\infty\setminus(\calA_\infty)_X}H$ and 
a generic line $L$ passing through $p$. 
Then the small loop $\gamma$ on $L$ around $p\in L$ determines 
a total turn monodromy $\rho(\gamma)\in GL_r(\bbC)$. 
The conjugacy class of $\rho(\gamma)$ in $GL_r(\bbC)$ 
depends only on $X\in L(\calA_\infty)$, which is denoted 
by $T_X$. 

The following vanishing theorem of local system cohomology groups 
is obtained in \cite{cdo}. 

\begin{Thm}
\normalfont
\label{thm:vanishing}
Let $\calL$ be a nonresonant local system on $M(\calA)$ of rank $r$, that is, 
for each dense subspace $X\subset H_\infty$ the corresponding 
monodromy operator $T_X$ does not admit $1$ as an eigenvalue. 
Then 
$$
\dim H^k(M(\calA), \calL)=
\left\{
\begin{array}{ll}
(-1)^\ell r\cdot\chi(M(\calA)) &\mbox{for }k=\ell \\
0 &\mbox{for }k\neq \ell , 
\end{array}
\right.
$$
where $\chi(M(\calA))$ is the Euler characteristic of 
the space $M(\calA)$. 
\end{Thm}
Note that $\calL$ is nonresonant if and only if the 
dual local system $\calL^\lor$ is nonresonant. 
From the universal coefficient theorem 
$$
H^k(M(\calA), \calL)\cong\Hom_\bbC(H_k(M(\calA), \calL^\lor), \bbC),
$$
we also have the similar vanishing theorem for 
local system homology groups $H_k(M(\calA), \calL)$. 

From Theorem \ref{thm:hamm} (a) 
we may identify, up to homotopy equivalence, 
$M(\calA)$ with a finite $\ell$-dimensional 
CW complex for which the 
\begin{equation}
\label{eq:1}
(\ell-1)\mbox{-skeleton has the homotopy type of } M(\calA)\cap U. 
\end{equation} 
We denote the attaching maps of $\ell$-cells by 
$\phi_k: \partial c_k\cong \rmS^{\ell-1}\rightarrow M(\calA)\cap U$, ($k=1, 
\cdots, b=b_\ell(M(\calA))$), 
where $c_k\cong\rmD^\ell$ is the $\ell$-dimensional unit disk. 
Hence $\phi=\{\phi_k\}_{k=1, \cdots, b}$ satisfies 
$$
\left(\left(M(\calA)\cap U\right)\cup_\phi\bigcup_k c_k\right)
\mbox{ is homotopic to }M(\calA). 
$$

Let $\calL$ be a rank $r$ local system over $M=M(\calA)$. 
For our purposes, 
it suffices to prove that $h(\phi_k)$
($k=1, \cdots, b$) generate $H_{\ell- 1}(M(\calA)\cap U, i^*\calL)$. 

Let 
\begin{equation}
\label{eq:chain1}
0\longrightarrow 
C_\ell 
\stackrel{\partial_\calL}{\longrightarrow}
C_{\ell-1}
\stackrel{\partial_\calL}{\longrightarrow}
\cdots
\stackrel{\partial_\calL}{\longrightarrow}
C_0
\longrightarrow 
0
\end{equation}
be the twisted cellular chain complex 
associated with the CW decomposition for $M(\calA)$. 
Then from (\ref{eq:1}), 
the twisted chain complex for $M(\calA)\cap U$ is obtained by 
truncating (\ref{eq:chain1}) as 
\begin{equation}
\label{eq:chain2}
0\longrightarrow 
C_{\ell-1}
\stackrel{\partial_\calL}{\longrightarrow}
\cdots
\stackrel{\partial_\calL}{\longrightarrow}
C_0
\longrightarrow 
0. 
\end{equation}
It is easily seen that if $\calL$ is generic in the sense of 
Theorem \ref{thm:vanishing}, then the restriction 
$i^*\calL$ is also generic. 
Applying Theorem \ref{thm:vanishing} to (\ref{eq:chain1}), only 
the $\ell$-th homology survives. Similarly, only 
the $(\ell-1)$-st homology survives in (\ref{eq:chain2}). 
Note that 
$H_{\ell-1}(M(\calA)\cap U, i^*\calL)=
{\mathrm{Ker}}(\partial_\calL: C_{\ell-1}\rightarrow C_{\ell-2})$. 
Thus we conclude that 
\begin{equation}
\label{eq:bd}
\partial_\calL:C_\ell\longrightarrow H_{\ell-1}(M(\calA)\cap U, 
i^*\calL)
\end{equation}
is surjective. 
Since the map (\ref{eq:bd}) is determined by 
$$
C_\ell\ni [c_k]\longmapsto [\partial c_k]=h(\phi_k), 
$$
$\{h(\phi_k)\}_{k=1, \cdots, b}$ 
generate $H_{\ell-1}(M(\calA)\cap U, i^*\calL)$. 
This completes the proof of Theorem \ref{thm:main}.

\begin{Lemma}
\label{lem:nonzero}
\normalfont
The Euler characteristic of $M(\calA)\cap U$ is not equal to zero, 
more precisely, 
$$
(-1)^{\ell-1}\chi(M(\calA)\cap U)>0. 
$$ 
\end{Lemma}

Given a hyperplane $H\in\calA$, we define $\calA'=\calA\setminus\{H\}$ 
and $\calA''=\calA'\cap H$. Then characteristic polynomials for 
these arrangements satisfy an inductive formula: 
$$
\chi(\calA, t)=\chi(\calA', t)-\chi(\calA'', t). 
$$

By Theorem \ref{thm:os}, 
the Euler characteristic $\chi(M(\calA))$ of the 
complement is equal to $\chi(\calA, 1)$.

\noindent
{\bf Proof of the Lemma \ref{lem:nonzero}.} 
From (\ref{eq:trunc}) and definition of the characteristic 
polynomial, we have 
$$
\chi(\calA\cap U, t)=\frac{\chi(\calA, t)-\chi(\calA, 0)}{t}. 
$$
The proof of the lemma is by induction on the number of hyperplanes. 
If $|\calA|=\ell$, $\calA$ is linearly isomorphic to the Boolean 
arrangement, i.e. one defined by 
$\{x_1\cdot x_2\cdots x_\ell=0\}$, 
for a certain coordinate system $(x_1, \cdots, x_\ell)$. 
In this case, $\chi(\calA, t)=(t-1)^\ell$, and we have 
$(-1)^{\ell-1}\chi(M(\calA)\cap U)=1$. Assume that $\calA$ contains 
more than $\ell$ hyperplanes. We can choose a hyperplane $H\in\calA$ 
such that $\calA'=\calA\setminus\{H\}$ is essential. Then 
$\calA''=\calA'\cap H$ is also essential, and obviously $U$ is 
generic to $\calA'$ and $\calA''$. Thus we have 
\begin{eqnarray*}
(-1)^{\ell-1}\chi(\calA\cap U)&=&(-1)^{\ell-1}\chi(\calA\cap U, 1)\\
&=&(-1)^{\ell-1}\left(\chi(\calA'\cap U, 1)-\chi(\calA''\cap U, 1)\right)\\
&=&(-1)^{\ell-1}\chi(\calA'\cap U, 1)+(-1)^{\ell-2}\chi(\calA''\cap U, 1)\\
&>&0. 
\end{eqnarray*}

\owari

Using Lemma \ref{lem:nonzero}, we have the following 
non-vanishing of the homotopy group, which generalizes 
a classical result of Hattori \cite{hat}. 

\begin{Cor}
\label{cor:homotopy}
\normalfont
Let $2\leq k\leq \ell-1$ and 
$F^k\subset V$ be a $k$-dimensional subspace generic to $\calA$. 
Then 
$
\pi_k(M(\calA)\cap F^k)\neq 0. 
$ 
\end{Cor}

\begin{Rem}
\label{rem:remark}
\normalfont
We can also prove Corollary \ref{cor:homotopy} 
directly in the following way. 
Suppose 
$
\pi_{\ell-1}(M(\calA)\cap U)=0. 
$ 
Then the attaching maps $\{\phi_k:\partial c_k=S^{\ell-1}\rightarrow M(\calA)\cap U\}$ 
of the top cells are homotopic to the constant map. 
Hence we have a homotopy equivalence 
$$
M(\calA)\mbox{ is homotopic to }
(M(\calA)\cap U)\vee \bigvee_k S^{\ell}. 
$$
However this contradicts to the fact that cohomology ring 
$H^*(M(\calA), \bbZ)$ is generated by degree one elements 
(Theorem \ref{thm:os}). Hence we have $
\pi_{\ell-1}(M(\calA)\cap U)\neq 0. 
$ 

\end{Rem}

\begin{Rem}
\normalfont
We should also note that other results on the non-vanishing 
of higher homotopy groups of generic sections are 
found in Randell \cite{ran-hom} (for generic sections of 
aspherical arrangements), in Papadima-Suciu \cite{pap-suc} 
(for hypersolvable arrangements) and in Dimca-Papadima 
\cite{dim-pap} (for iterated generic hyperplane sections). 
\end{Rem}

\vspace{4mm}

\noindent
{\bf Acknowledgements. }
The author expresses deep gratitude to 
Professor M. Falk and Professor S. Papadima, 
for indicating Remark \ref{rem:remark} to him. 
The author also thanks to the referee for 
a lot of suggestions which improve the paper. 
The author was supported by the JSPS Postdoctoral 
Fellowships for Research Abroad.

\end{document}